\documentclass[10pt]{amsart}
\usepackage{amsfonts,amssymb,stmaryrd,amscd,amsmath,latexsym,amsbsy}

\usepackage{amssymb}
\usepackage{amsfonts}
\usepackage{latexsym}
\usepackage{pdfsync}
\usepackage{color}
\newtheorem{theorem}{Theorem}[section]
\newtheorem{lemma}[theorem]{Lemma}
\newtheorem{proposition}[theorem]{Proposition}

\theoremstyle{definition}
\newtheorem{definition}[theorem]{Definition}

\newtheorem{remark}[theorem]{Remark}

\newcommand{\End}{\text{End}}
\newcommand{\h}{\mathfrak{h}}
\newcommand{\Hom}{\text{Hom}}
\newcommand{\tr}{\text{tr}\,}
\newcommand{\kk}{\mathfrak{k}}
\def\HH{\hbox{${\mathcal H}$\kern-5.2pt${\mathcal H}$}}

\begin{document}
\title[Construction of representations of $d$AHA and $d$DAHA of 
type $BC_n$]{A Lie-theoretic
construction of some representations of the degenerate affine and
double affine Hecke algebras of type $BC_n$}
%%%%%%%%%%%%%%%%%%%%%%%%

\author[Pavel Etingof]{Pavel Etingof}
\author[Rebecca Freund]{Rebecca Freund}
\author[Xiaoguang Ma]{Xiaoguang Ma}

\address{Department of Mathematics, Massachusetts Institute of Technology,
Cambridge, MA 02139, USA}
\email{etingof@math.mit.edu}

\address{Massachusetts Institute of Technology,
Cambridge, MA 02139, USA}
\email{rlfreund@mit.edu}

\address{Department of Mathematics, Massachusetts Institute of Technology,
Cambridge, MA 02139, USA}
\email{xma@math.mit.edu}

%%%%%%%%%%%%%

\maketitle

\section{Introduction}
%%%%%%%%%%%%%%%

The degenerate affine Hecke algebra (dAHA) of any Coxeter group
was defined by Drinfeld and Lusztig(\cite{Dri},\cite{Lus}). 
It is generated by the group algebra of the Coxeter group and by the
commuting generators $y_{i}$ with some relations.

The degenerate double affine Hecke algebra (dDAHA) of a root system 
was introduced by Cherednik (see \cite{Ch}). It is generated by the the 
group algebra of the Weyl group, commuting generators $y_{i}$, and 
by another kind of commuting generators $X_{i}$ with some relations.
The dAHA corresponding to the Weyl group can be realized 
as a subalgebra of dDAHA, generated by the Weyl group and the
elements $y_i$. 

The paper \cite{AS} gives a Lie-theoretic construction of
representations of the dAHA of type $A_{n-1}$. Namely (see
\cite{CEE}, section 9), for every $\mathfrak{sl}_N$-bimodule
$M$, an action of the dAHA $\mathcal{H}$ of type $A_{n-1}$ is
constructed on the space
\begin{equation*}
F_n(M):=(M\otimes (\mathbb{C}^{N})^{\otimes n})^{\mathfrak{sl}_N}, 
\end{equation*}
where the invariants are taken with respect to the adjoint action
of $\mathfrak{sl}_N$ on $M$. 

This construction is upgraded to a Lie-theoretic construction of
representations of dDAHA of type $A_{n-1}$ 
\cite{CEE}, Section 9. Namely, for any 
$\mathcal{D}$-module $M$ on $SL_{N}$, 
the paper \cite{CEE} constructs an action of dDAHA $\HH$ with parameter
$k=N/n$ on the space $F_n(M)$, such that the induced action of
the dAHA $\mathcal{H}\subset \HH$
coincides with the action of \cite{AS}, obtained by regarding $M$ as an
$\mathfrak{sl}_N$-bimodule via left- and right-invariant vector
fields on $SL_{N}$.

The main result of this paper is an analog of the constructions
of \cite{AS} and \cite{CEE} for dAHA and dDAHA of type $BC_n$,
which gives a method of obtaining representations of these
algebras from Lie theory. Specifically, given a
module $M$ over the Lie algebra 
$\mathfrak{g}:=\mathfrak{gl}_N$, we first construct an action of
the dAHA $\mathcal{H}$ of type $B_{n}$ on the space
$F_{n,p,\mu}(M)$ of $\mu$-invariants in $M\otimes
(\mathbb{C}^{N})^{\otimes n}$ under the subalgebra
$\kk_0:=(\mathfrak{gl}_p\oplus \mathfrak{gl}_q)\cap 
\mathfrak{sl}_N\subset \mathfrak{g}$, where $q=N-p$, and
$\mu\in \Bbb C$ is a parameter (here by $\mu$-invariants we mean 
eigenvectors of $\kk_0$ with eigenvalues given by 
the character $\mu\chi$, where $\chi$ is a basic character of
$\kk_0$). In this construction, the
parameters of $\mathcal{H}$ are certain explicit functions of
$\mu$ and $p$. Thus we obtain an functor $F_{n,p,\mu}$ from the
category of $\mathfrak{gl}_N$-modules to the category of
representations of $\mathcal{H}$. It is easy to see that this
functor factors through the category of Harish-Chandra modules
for the symmetric pair $(\mathfrak{gl}_N,\mathfrak{gl}_p\oplus
\mathfrak{gl}_q)$, so it suffices to restrict our attention to
Harish-Chandra modules.
In particular, 
the resulting functor after this restriction, i.e. the functor
from the category of Harish-Chandra modules to the category of
representations of $\mathcal{H}$, is exact.
\footnote{The functor $F_{n,p,\mu}$ is not exact without this restriction. 
The reason is that
the functor of $\mathfrak{g}$-invariants for a semisimple Lie algebra $\mathfrak{g}$ is only exact
on the category of locally finite $\mathfrak{g}$-modules.}

Then we upgrade this construction to one giving representations
of dDAHA $\HH$ of type $BC_n$. Namely, let $G=GL_N$, and
$K=GL_p\times GL_q\subset G$. Then for any
${\lambda}$-twisted $\mathcal{D}$-module $M$ on $G/K$ we
construct an action of the dDAHA $\HH$ of type $BC_n$ on
the space $F_{n,p,\mu}(M)$. In this construction, the parameters
of $\HH$ are certain explicit functions of $\lambda$, $\mu$, and
$p$. Moreover, the underlying representation of $\mathcal{H}$
coincides with the representation obtained in the previous
construction, if we regard $M$ as a $\mathfrak{gl}_N$-module
via the vector fields corresponding to the action of $G$ on
$G/K$.  Thus we obtain an functor $F^\lambda_{n,p,\mu}$ from the
category of $\lambda$-twisted $\mathcal{D}$-modules on $G/K$ to
the category of representations of $\HH$. This functor factors
through the category of $K$-monodromic twisted
$\mathcal{D}$-modules, so it suffices to restrict our attention
to such $\mathcal{D}$-modules. In particular, 
the resulting functor after this restriction, i.e. the functor
from the category of $K$-monodromic twisted
$\mathcal{D}$-modules to the category of representations of $\HH$, is exact.
\footnote{Similarly to the affine case, $F^\lambda_{n,p,\mu}$ is not exact without this restriction.}

Since the appearance of the first version of this paper,
the functor $F_{n,p,\mu}$ has been studied in the followup paper
\cite{M} by the third author. In this paper, it was shown that 
the principal series representations of $U(p,q)$ are mapped by 
the $F_{n,p,\mu}$ to certain induced modules over the dAHA. 
This is analogous to the result of \cite{AS}, where it is shown that 
in the type $A$ case, standard modules go to standard modules. 

We expect that further careful study of the functors $F_{n,p,\mu}$ and 
$F_{n,p,\mu}^\lambda$ will reveal an interesting new connection
between, on the one hand, the representation theory of  
(the universal cover of) the Lie group $U(p,q)$ and the theory of
monodromic $\mathcal{D}$-modules on the (complexification of the) 
corresponding symmetric space $U(p,q)/U(p)\times U(q)$, and, on
the other hand, the representation theory of dAHA and dDAHA of 
type $BC_n$. We plan to discuss this in future publications. 

The paper is organized as follows.  In Section 2, we recall the
definitions of dAHA and dDAHA of a root system. In Section 3, 
we write down the explicit definitions for dAHA and dDAHA of type
$BC_{n}$. In Section 4, we construct the
functor $F_{n,p,\mu}$. In Section 5, we construct the functor 
$F_{n,p,\mu}^\lambda$. In Section 6, we study some properties of
this functor. 
 
\section{Definitions and notations}
%%%%%%%%%%%%%%%%%%%

In this section, we recall the definitions 
of degenerate affine and double affine Hecke algebras. 
For more details, see \cite{Ch}.

Let $\mathfrak{h}$ be a finite dimensional real vector space with a
positive definite symmetric bilinear form $(\cdot,\cdot)$.  Let
$\{\epsilon_{i}\}$ be a basis for $\mathfrak{h}$ such that
$(\epsilon_{i},\epsilon_{j})=\delta_{ij}$.  Let $R$ be an irreducible root
system in $\mathfrak{h}$ (possibly non-reduced). 
Let $R_{+}$ be the set of positive roots of $R$, and
let $\Pi=\{\alpha_{i}\}$ be the set of simple roots.

For any root $\alpha$, the corresponding 
coroot is $\alpha^{\vee}=2\alpha/(\alpha,\alpha)$. 
Let $Q$ and $Q^{\vee}$ be the root lattice and the coroot lattice. Let 
$P=\Hom_{\mathbb{Z}}(Q^{\vee},\mathbb{Z})$ be the weight lattice. 
 
Let $W$ be the Weyl group of $R$. 
Let $\Sigma$ be the set of reflections in $W$. Let 
 $S_\alpha\in \Sigma$ be the reflection corresponding to the root
$\alpha$. In particular, we write $S_i$ for the simple
reflections $S_{\alpha_i}$.

If $\alpha$ is a root, let $\nu_\alpha=1$ if $\alpha$ an indivisible root, 
and $\nu_\alpha=2$ otherwise. 

Let us define dAHA. 
Let $\kappa: \Sigma\to \Bbb C$ be a conjugation invariant
function. 

\begin{definition}
The {\em degenerate affine Hecke algebra (dAHA)}
$\mathcal{H}(\kappa)$ is the quotient of the free product 
$\Bbb CW * S\mathfrak{h}$ by the 
relations
$$
S_i y-y^{S_i} S_i=\kappa(S_i)\alpha_i(y),\quad y\in \mathfrak{h}.
$$
\end{definition}
For any $a\ne 0$, multiplication of $y$-generators by $a$ defines an isomorphism of $\mathcal{H}(\kappa)$ with $\mathcal{H}(a\kappa)$
(under this transformation, nothing happens to $\alpha_{i}(y)$ since they are just numbers).
Thus in the simply laced situation, there is only one 
nontrivial case, $\kappa=1$, and in the non-simply laced case, 
the function $\kappa$ takes 
two values $\kappa_1,\kappa_2$ (the values of $\kappa$ on
the root reflections for long and short indivisible roots, 
respectively), and the algebra depends only
on the ratio $\kappa_2/\kappa_1$ (unless both values are zero). 

Now let us define dDAHA.  
Let $k: R\to \Bbb C$, $\alpha\mapsto k_\alpha$, be a 
function such that $k_{g(\alpha)}=k_{\alpha}$ for all $g\in
W$. Let $t\in \Bbb C$. 

\begin{definition}
For $\epsilon\in \mathfrak{h}$, 
define the Dunkl-Cherednik operator 
\begin{equation*}
D_{\epsilon}(t,k)=t\partial_{\epsilon}-\sum_{\alpha\in R_{+}}
\frac{k_{\alpha}\alpha(\epsilon)}{1-e^{-\alpha}}(1-S_{\alpha})
+\rho(k)(\epsilon),
\end{equation*}
where $\partial_{\epsilon}$ is the differentiation along
$\epsilon$, and $\rho(k)=\frac{1}{2}\sum_{\alpha\in R_{+}}k_{\alpha}\alpha$. 
This operator acts on the space $E$ 
of trigonometric polynomials on $\h/Q^\vee$.  
\end{definition}

An important property of the operators $D_{\epsilon}$ is that
they commute with each other.

\begin{definition}
The {\em degenerate double affine Hecke algebra (dDAHA)} 
$\HH(t,k)$ is generated by $W$, the Dunkl-Cherednik operators, 
and the elements $e^{\lambda}\quad(\lambda\in P)$.
\end{definition}

\begin{remark}
This is not the original definition of dDAHA. But it is
equivalent to the original definition by a theorem of Cherednik,
see \cite{Ch}. 
\end{remark}

Obviously, for any $a\ne 0$, we have a natural isomorphism 
between $\HH(t,k)$ and $\mathcal{H}(at,ak)$;
thus, there are only two essentially different cases:
$t=0$ (the classical case) and $t=1$ (the quantum case).  

The following proposition can be proved by a straightforward
computation. 

\begin{proposition} (see \cite{Ch})
The subalgebra of the dDAHA $\HH(t,k)$ generated by $W$ and the
Dunkl-Cherednik operators is isomorphic to the dAHA 
$\mathcal{H}(\kappa)$, where $\kappa(S)=\sum_{\alpha:
S=S_\alpha}k_\alpha\nu_\alpha$.
\end{proposition}

\section{Type $BC_n$ $d$AHA and $d$DAHA}
%%%%%%%%%%%%%%%%%

\subsection{Definitions of the type $BC_n$ dAHA and dDAHA}
%%%%%%%%%%%%%%%%%%%%%%%%%%%%%%

Let us now describe the dAHA and dDAHA of type $BC_n$ more
explicitly.

Let $\mathfrak{h}$ be a real vector space of dimension $n$ with
orthonormal basis $\epsilon_1, \ldots, \epsilon_n$. 
We will identify $\mathfrak{h}$ with its dual by the bilinear
form, and set $X_{i}=e^{\epsilon_{i}}$, $y_{i}=D_{\epsilon_{i}}$.

The roots of
type $BC_{n}$ are
\begin{equation*}
R=\{\pm\epsilon_{i}\}\cup\{\pm2\epsilon_{i}\}
\cup\{\pm\epsilon_{i}\pm\epsilon_{j}\}_{i\neq j},
\end{equation*}
and the positive roots are 
\begin{equation*}
R_{+}=\{\epsilon_{i}\}\cup\{2\epsilon_{i}\}\cup\{\epsilon_{i}\pm\epsilon_{j}\}_{i< j}.
\end{equation*}

The function $\kappa$ considered in the previous section
reduces to two parameters $\kappa=(\kappa_1,\kappa_2)$,  
while the function $k$  reduces
to three parameters $k=(k_{1},k_{2},k_{3})$ 
corresponding to the three kinds of positive roots: those of  
lengths $2,1,4$, respectively. Namely, $k_1=k_{\epsilon_i-\epsilon_j}$,
$k_2=k_{\epsilon_i}$, $k_3=2k_{2\epsilon_i}$.  

Let $W=S_n\ltimes (\Bbb Z/2\Bbb Z)^n$ 
be the Weyl group of type $BC_{n}$.
We denote by $S_{ij}$ the reflection in this group 
corresponding to the root $\epsilon_i-\epsilon_j$, 
and by $\gamma_{i}$ the reflection corresponding to
$\epsilon_{i}$. Then $W$ is generated by 
$S_{i}=S_{i,i+1}, i=1,\ldots, n-1 $ and $\gamma_{n}$. 

The type $BC_{n}$ 
dAHA $\mathcal{H}(\kappa_1,\kappa_2)$
is then defined as follows: 
\begin{itemize}
\item generators: $y_{1},\ldots,y_{n}$ 
and $\mathbb{C}W$;

\item relations: 
\begin{enumerate}
\item[i)] $S_{i}$ and $\gamma_{n}$ satisfy the Coxeter relations;
\item[ii)] $S_{i}y_{i}-y_{i+1}S_{i}=\kappa_{1}$,
$[S_{i},y_{j}]=0, (j\neq i,i+1)$;
\item[iii)]  $\gamma_{n}y_{n}+y_{n}\gamma_{n}=\kappa_2$, 
$[\gamma_{n},y_{j}]=0, (j\neq n)$;
\item[iv)] $[y_{i},y_{j}]=0$. 
\end{enumerate}
\end{itemize}

On the other hand, the type $BC_n$ 
dDAHA $\HH(t,k_{1},k_{2},k_{3})$ is
defined as follows:

\begin{itemize}
\item generators: $X_{1}, \ldots, X_{n}$, $y_{1},\ldots,y_{n}$ 
and $\mathbb{C}W$;

\item relations: 
\begin{enumerate}
\item[i)] $S_{i}$ and $\gamma_{n}$ satisfy the Coxeter relations;
\item[ii)] $S_{i}X_{i}-X_{i+1}S_{i}=0$, $[S_{i},X_{j}]=0, (j\neq i,i+1)$;
\item[iii)] $S_{i}y_{i}-y_{i+1}S_{i}=k_{1}$, $[S_{i},y_{j}]=0, (j\neq i,i+1)$;
\item[iv)]  $\gamma_{n}y_{n}+y_{n}\gamma_{n}=k_{2}+k_{3}$, 
$\gamma_{n}X_{n}=X_{n}^{-1}\gamma_{n}$, 
\\$[\gamma_{n},y_{j}]=[\gamma_{n},X_{j}]=0, (j\neq n)$;
\item[v)] $[X_{i},X_{j}]=[y_{i},y_{j}]=0$; 
\item[vi)]$[y_{j},X_{i}]=k_{1}X_{i}S_{ij}-k_{1}X_{i}S_{ij}\gamma_{i}
\gamma_{j}$, \\ $[y_{i},X_{j}]=k_{1}X_{i}S_{ij}-k_{1}X_{j}S_{ij}\gamma_{i}\gamma_{j},(i<j)$;
\item[vii)] 
\begin{eqnarray*}
[y_{i},X_{i}]
&=&tX_{i}-k_{1}X_{i}\sum_{k>i}S_{ik}-k_{1}\sum_{k<i}S_{ik}X_{i}-k_{1}X_{i}\sum_{k\neq i}S_{ik}\gamma_{i}\gamma_{k}\\
&&\qquad- (k_{2}+k_{3})X_{i}\gamma_{i}-k_{2}\gamma_{i}.
\end{eqnarray*}
\end{enumerate}
\end{itemize}

In particular, we see that the subalgebra in the dDAHA generated 
by $W$ and $y_i$ is $\mathcal{H}(\kappa_1,\kappa_2)$, where 
$\kappa_1=k_1$ and $\kappa_2=k_2+k_3$.

\subsection{Another set of generators of dAHA of type $BC_{n}$}

The advantage of the generators $y_i$ is that they commute with
each other, but their disadvantage is that they do not change
according to the standard representation of the Weyl group. 
It turns out that it is possible (and useful in some situations, 
including one of this paper) to trade the first property for the
second one, by replacing the generators $y_i$ by their shifted
versions $\tilde y_i$, passing from Lusztig's presentation 
of dAHA (\cite{Lus}) to Drinfeld's one (\cite{Dri}).
\footnote{See \cite{RS} for more details.}  

Namely, for each $i=1,\ldots, n$, define 
\begin{equation*}
\tilde{y}_{i}=y_{i}-\frac{\kappa_{2}}{2}\gamma_{i}-
\frac{\kappa_{1}}{2}\sum_{k>i}S_{ik}+\frac{\kappa_{1}}{2}\sum_{k<i}S_{ik}-
\frac{\kappa_{1}}{2}\sum_{i\neq k}S_{ik}\gamma_{i}\gamma_{k}.
\end{equation*}

\begin{lemma}\label{rel-tildey1}
The type $BC_{n}$ dAHA $\mathcal{H}(\kappa_{1},\kappa_{2})$ is generated by
$w\in W$ and $\tilde{y}_{i}$ with the following relations:
\begin{enumerate}
\item[i)] $S_{i}\tilde{y}_{i}-\tilde{y}_{i+1}S_{i}=0$, 
\quad $S_{j}\tilde{y}_{i}-\tilde{y}_{i}S_{j}=0$, $(i\neq j)$;
\item[ii)] $\tilde{y}_{n}\gamma_{n}+\gamma_{n}\tilde{y}_{n}=0$, 
\quad $\tilde{y}_{i}\gamma_{n}-\gamma_{n}\tilde{y}_{i}=0$, $(i\neq n)$;
\item[vi)]
\begin{eqnarray*}\label{ycom1}
[\tilde{y}_i, \tilde{y}_j]
&=&  \frac{\kappa_{1}\kappa_{2}}{2}S_{ij}(\gamma_{j}
-\gamma_{i})
  + \frac{\kappa_{1}^{2}}{4}\sum_{k \neq i, j}S_{jk}S_{ik} 
  - \frac{\kappa_{1}^{2}}{4}\sum_{k \neq i, j}S_{ik}S_{jk}  \\ 
&& \quad + \frac{\kappa_{1}^{2}}{4}
\sum_{k \neq i, j}S_{ik}S_{jk}
(-\gamma_i\gamma_k + \gamma_i\gamma_j+\gamma_j\gamma_k)\\
&& 
\qquad- \frac{\kappa_{1}^{2}}{4}\sum_{k \neq i, j}
S_{jk}S_{ik}(\gamma_i\gamma_j-\gamma_j\gamma_k+\gamma_i\gamma_k). 
\end{eqnarray*}
\end{enumerate}
\end{lemma}

\begin{proof}
The proof is contained in the proof of Lemma \ref{rel-tildey}
given in the next subsection. 
\end{proof}

\subsection{Another set of generators of dDAHA of type $BC_{n}$}

Similarly to the previous subsection, 
for each $i=1,\ldots, n$, define 
\begin{equation*}
\tilde{y}_{i}=y_{i}-\frac{k_{2}+k_{3}}{2}\gamma_{i}-
\frac{k_{1}}{2}\sum_{k>i}S_{ik}+\frac{k_{1}}{2}\sum_{k<i}S_{ik}-
\frac{k_{1}}{2}\sum_{i\neq k}S_{ik}\gamma_{i}\gamma_{k}
\end{equation*}
(together with the Weyl group, the elements $\tilde y_i$
generate dAHA in Drinfeld's presentation \cite{Dri}).

\begin{lemma}\label{rel-tildey}
The type $BC_{n}$ $d$DAHA $\HH(t,k_{1},k_{2},k_{3})$ is generated by
$w\in W$, $X_{i}$ and $\tilde{y}_{i}$ with the following relations:
\begin{enumerate}
\item[i)] the relations among elements of $W$ and $X_{i}$ are the same as before;
\item[ii)] $S_{i}\tilde{y}_{i}-\tilde{y}_{i+1}S_{i}=0$, 
\quad $S_{j}\tilde{y}_{i}-\tilde{y}_{i}S_{j}=0$, $(i\neq j)$;
\item[iii)] $\tilde{y}_{n}\gamma_{n}+\gamma_{n}\tilde{y}_{n}=0$, 
\quad $\tilde{y}_{i}\gamma_{n}-\gamma_{n}\tilde{y}_{i}=0$, $(i\neq n)$;
\item[iv)] $[\tilde{y}_{j},X_{i}]=\dfrac{k_{1}}{2}(X_{i}+X_{j})S_{ij}-
\dfrac{k_{1}}{2}(X_{i}+X_{j}^{-1})S_{ij}\gamma_{i}\gamma_{j}$,
$(i\ne j)$;
\item[v)] 
\begin{eqnarray*}
[\tilde{y}_{i},X_{i}] 
& = & tX_{i}-\frac{k_{2}+k_{3}}{2}(X_{i}^{-1}+X_{i})\gamma_{i}-
k_{2}\gamma_{i}\\
&&\qquad-\frac{k_{1}}{2}\sum_{k\neq i}((X_{i}+X_{k})S_{ik}+(X_{i}+X_{k}^{-1})S_{ik}\gamma_{i}\gamma_{k});
\end{eqnarray*}
\item[vi)]
\begin{eqnarray*}\label{ycom}
[\tilde{y}_i, \tilde{y}_j]
&=&  \frac{k_{1}(k_{2}+k_{3})}{2}S_{ij}(\gamma_{j}
-\gamma_{i})
  + \frac{k_{1}^{2}}{4}\sum_{k \neq i, j}S_{jk}S_{ik} 
  - \frac{k_{1}^{2}}{4}\sum_{k \neq i, j}S_{ik}S_{jk}  \\ 
&& \quad + \frac{k_{1}^{2}}{4}\sum_{k \neq i, j}S_{ik}S_{jk}(-\gamma_i\gamma_k + \gamma_i\gamma_j+\gamma_j\gamma_k)\\
&& 
\qquad- \frac{k_{1}^{2}}{4}\sum_{k \neq i, j}S_{jk}S_{ik}(\gamma_i\gamma_j-\gamma_j\gamma_k+\gamma_i\gamma_k), (i \ne j). 
\end{eqnarray*}
\end{enumerate}

\end{lemma}

\begin{proof}
Only the last relation is nontrivial. Its proof is by a direct
computation. To make formulas more compact, 
set 
$$
R_{i}=-\frac{k_{1}}{2}\sum_{k>i}S_{ik}
+\frac{k_{1}}{2}\sum_{k<i}S_{ik}
-\frac{k_{1}}{2}\sum_{i\neq k}S_{ik}\gamma_{i}\gamma_{k},
$$ 
then we have
$$
y_{i}=\tilde{y}_{i}+\frac{k_{2}+k_{3}}{2}\gamma_{i}
-R_{i}.
$$

Assume $i<j$, then we have
\begin{eqnarray*} 
&&[y_{i},y_{j}]\\
&&\quad=[\tilde{y}_{i},\tilde{y}_{j}]-[\tilde{y}_{i},R_{j}]-[R_{i},\tilde{y}_{j}]-\frac{k_{2}+k_{3}}{2}([\gamma_{i},R_{j}]+[R_{i},\gamma_{j}])
+[R_{i},R_{j}].
\end{eqnarray*} 

Since
\begin{eqnarray*}
&&[\tilde{y}_{i},R_{j}]
=\frac{k_{1}}{2}S_{ij}\tilde{y}_{j}
-\frac{k_{1}}{2}S_{ij}\tilde{y}_{i}
-\frac{k_{1}}{2}S_{ij}\tilde{y}_{j}\gamma_{i}\gamma_{j}
-\frac{k_{1}}{2}S_{ij}\tilde{y}_{i}\gamma_{i}\gamma_{j},\\
&&[R_{i},\tilde{y}_{j}]
=\frac{k_{1}}{2}S_{ij}\tilde{y}_{i}
-\frac{k_{1}}{2}S_{ij}\tilde{y}_{j}
+\frac{k_{1}}{2}S_{ij}\tilde{y}_{j}\gamma_{i}\gamma_{j}
+\frac{k_{1}}{2}S_{ij}\tilde{y}_{i}\gamma_{i}\gamma_{j},\\
&&[\gamma_{i},R_{j}]
=k_{1}S_{ij}\gamma_{j}
-k_{1}S_{ij}\gamma_{i},\\
&&[R_{i},\gamma_{j}]
=0,
\end{eqnarray*}
we have
$$[\tilde{y}_{i},\tilde{y}_{j}]
=\frac{k_{1}(k_{2}+k_{3})}{2}S_{ij}(\gamma_{j}
-\gamma_{i})
-[R_{i},R_{j}].
$$
By direct computation, we have $[R_i,R_j]=\frac{k_1^2}{4}R_{ij}$, where 
\begin{eqnarray*}
R_{ij}
&=&  \sum_{k>j}S_{ij}S_{jk} 
- \sum_{k>j}S_{jk}S_{ij}
+  \sum_{k>j} S_{ik}S_{jk} 
- \sum_{k>j}S_{jk}S_{ik} \\
&&-\mathop{\sum_{i<k<j}}_{\text{or }k>j}S_{ik}S_{ij} 
+  \mathop{\sum_{i<k<j}}_{ \text{or }k>j}S_{ij}S_{ik}
- \mathop{\sum_{k<i\text{ or}}}_{i<k<j}S_{ij}S_{jk} 
+ \mathop{\sum_{k<i\text{ or}}}_{i<k<j}S_{jk}S_{ij}\\ 
&&-\sum_{i<k<j}S_{ik}S_{jk} 
+ \sum_{i<k<j}S_{jk}S_{ik}
+ \mathop{\sum_{i<k<j}}_{\text{or }k>j}S_{ik}S_{ij}\gamma_i\gamma_j\\
&&- \mathop{\sum_{i<k<j}}_{\text{or }k>j}S_{ij}S_{ik}\gamma_j\gamma_k
 +\mathop{\mathop{\sum_{k<i\text{ or}}}_{i<k<j}}_{\text{or }k>j}S_{ij}S_{jk}\gamma_j\gamma_k 
- \mathop{\mathop{\sum_{k<i \text{ or}}}_{i<k<j}}_{\text{or }k>j} S_{jk}S_{ij}\gamma_i\gamma_k\\ 
&&+\mathop{\sum_{i<k<j}}_{\text{or }k>j}S_{ik}S_{jk}\gamma_j\gamma_k 
- \mathop{\sum_{i<k<j}}_{\text{or }k>j}S_{jk}S_{ik}\gamma_i\gamma_j\\
&&+ \sum_{k<i}S_{ik}S_{ij} 
- \sum_{k<i}S_{ij}S_{ik} 
+ \sum_{k<i} S_{ik}S_{jk} 
- \sum_{k<i}S_{jk}S_{ik}\\ 
&& -\sum_{k<i}S_{ik}S_{ij}\gamma_i\gamma_j 
+ \sum_{k<i}S_{ij}S_{ik}\gamma_j\gamma_k 
- \sum_{k<i}S_{ik}S_{jk}\gamma_j\gamma_k 
+ \sum_{k<i} S_{jk}S_{ik}\gamma_i\gamma_j\\  
&&+ \sum_{k>j}S_{ij}S_{jk}\gamma_i\gamma_k 
- \sum_{k>j}S_{jk}S_{ij}\gamma_i\gamma_j 
+ \sum_{k>j}S_{ik}S_{jk}\gamma_i\gamma_j 
- \sum_{k>j}S_{jk}S_{ik}\gamma_i\gamma_k \\ 
&&- \mathop{\mathop{\sum_{k<i \text{ or}}}_{i<k<j}}_{\text{or }k>j}S_{ik}S_{ij}\gamma_j\gamma_k 
+ \mathop{\mathop{\sum_{k<i \text{ or}}}_{i<k<j}}_{\text{or }k>j}S_{ij}S_{ik}\gamma_i\gamma_k
- \mathop{\sum_{k<i \text{ or}}}_{i<k<j}S_{ij}S_{jk}\gamma_i\gamma_k \\
&&+ \mathop{\sum_{k<i \text{ or}}}_{i<k<j}S_{jk}S_{ij}\gamma_i\gamma_j 
-\mathop{\sum_{k<i \text{ or}}}_{i<k<j}S_{ik}S_{jk}\gamma_i\gamma_j 
+ \mathop{\sum_{k<i \text{ or}}}_{i<k<j}S_{jk}S_{ik}\gamma_i\gamma_k \\ 
&&+ \sum_{k \neq i, j}S_{ik}S_{ij}\gamma_i\gamma_k 
- \sum_{k \neq i, j}S_{ij}S_{ik}\gamma_i\gamma_j
+ \sum_{k \neq i, j}S_{ij}S_{jk}\gamma_i\gamma_j \\
&&- \sum_{k \neq i, j}S_{jk}S_{ij}\gamma_j\gamma_k 
+ \sum_{k \neq i, j} S_{ik}S_{jk}\gamma_i\gamma_k 
- \sum_{k \neq i, j}S_{jk}S_{ik}\gamma_j\gamma_k
\end{eqnarray*}
\begin{eqnarray*}
&=&-\sum_{k \neq i, j}S_{jk}S_{ik} 
+ \sum_{k \neq i, j}S_{ik}S_{jk}  \\ 
&&  - \sum_{k \neq i, j}S_{ik}S_{jk}(-\gamma_i\gamma_k - \gamma_i\gamma_j+\gamma_j\gamma_k) 
+ \sum_{k \neq i, j}S_{jk}S_{ik}(\gamma_i\gamma_j-\gamma_j\gamma_k+\gamma_i\gamma_k). 
\end{eqnarray*}

Thus, from the commutativity of $y_{i}$ and the above, we get
\begin{eqnarray*}
[\tilde{y}_i, \tilde{y}_j]
&=&  \frac{k_{1}(k_{2}+k_{3})}{2}S_{ij}(\gamma_{j}
-\gamma_{i})
  + \frac{k_{1}^{2}}{4}\sum_{k \neq i, j}S_{jk}S_{ik} 
  - \frac{k_{1}^{2}}{4}\sum_{k \neq i, j}S_{ik}S_{jk}  \\ 
&& \quad + \frac{k_{1}^{2}}{4}\sum_{k \neq i, j}S_{ik}S_{jk}(-\gamma_i\gamma_k + \gamma_i\gamma_j+\gamma_j\gamma_k)\\
&& 
\qquad- \frac{k_{1}^{2}}{4}\sum_{k \neq i, j}S_{jk}S_{ik}(\gamma_i\gamma_j-\gamma_j\gamma_k+\gamma_i\gamma_k), 
\end{eqnarray*}
as desired.

\end{proof}

\section{Construction of the functor $F_{n,p,\mu}$}
%%%%%%%%%%%%%%%%%%%%%%%%%%%%%%%%%%

\subsection{Notations}
%%%%%%%%%%%%%

We will use the following notations. Let $p,q\in \mathbb{N}$ 
and let $N=p+q$.
Let $E_{ij}$ be the $N$ by $N$ matrix which has a $1$ at the 
$(i,j)$-th position and 0 elsewhere.  Set
$J=\left(\begin{array}{cc}I_p &  \\ & -I_q\end{array}\right)$,
where $I_p$ is the identity matrix of size $p$. 

Let $\mathfrak{g}=\mathfrak{gl}_N(\mathbb{C})$ 
be the Lie algebra of $G=GL_N(\mathbb{C})$. Let 
$\mathfrak{k}=\mathfrak{gl}_p(\Bbb C)\times\mathfrak{gl}_q(\Bbb
C)$ be the Lie algebra of  $K=GL_p(\mathbb{C})\times GL_q(\mathbb{C})$. Let 
$\mathfrak{k}_{0}$ be the subalgebra of trace zero 
matrices in $\mathfrak{k}$.

Define a character $\chi$ of $\mathfrak{k}$ by
$$\chi(\left(\begin{array}{cc}X_1 & 0 \\0 &
X_2\end{array}\right))=
q\tr X_{1}-p\tr X_{2}.$$
For this $\chi$, we have the following obvious lemma.

\begin{lemma}\label{chi-mu}
We have $\chi(E_{ij})=0$ for $i\neq j$, and
$$
\chi (E_{ii})=\left\{\begin{array}{cc}q,  & i\leq p; \\
&\\
-p ,& i> p.\end{array}\right.
$$
In particular, $\chi(I_N)=0$.
\end{lemma}

\begin{remark}
The property $\chi(I_N)=0$ is very important for the future 
discussion. In fact, this is also the reason why we choose such $\chi$.
\end{remark}

We will also use the following summation notations:
$$
\sum_{i,\ldots,j}=\sum_{i=1}^{n}\cdots \sum_{j=1}^{n},\quad
 \sum_{i\ldots j}=\sum_{i=1}^{p}\cdots
\sum_{j=1}^{p}+\sum_{i=p+1}^{n}\cdots\sum_{j=p+1}^{n},
$$
 $$
\sum_{i\ldots j|k\ldots
l}=\sum_{i=1}^{p}\cdots\sum_{j=1}^{p}\sum_{k=p+1}^{n}\cdots
\sum_{l=p+1}^{n}+\sum_{k=1}^{p}\cdots
\sum_{l=1}^{p}\sum_{i=p+1}^{n}\cdots\sum_{j=p+1}^{n}.
$$ 
Thus we have two ranges of summation ($[1,p]$ and $[p+1,N]$), and 
the indices not separated by anything must be in the same range, 
while indices separated by a vertical line must be in 
different ranges. Indices separated by a comma are independent. 

\subsection{Construction of the functor $F_{n,p,\mu}$}
%%%%%%%%%%%%%%%%%%%%%%%%%%%%%

Let $Y$ be a $\mathfrak{k}_0$-module.
For any $\mu\in \mathbb{C}$,  
define the space of $\mu$-invariants 
$Y^{\mathfrak{k}_0,\mu}$ to be the space 
of those $v\in Y$ for which $xv=\mu\chi(x)v$ for all $x\in
\mathfrak{k}_0$. 

Let $V=\mathbb{C}^{N}$ be the vector representation of $\mathfrak{g}$.
Let $M$ be a $\mathfrak{g}$-module. Define 
\begin{equation*}
F_{n,p,\mu}(M)=(M\otimes V^{\otimes n})^{\mathfrak{k}_0,\mu}.
\end{equation*}

The Weyl group $W$ acts on $M\otimes V^{\otimes n}$
in the following way: the element $S_{ij}$ acts by exchanging 
the $i$-th and $j$-th factors, 
and $\gamma_{i}$ acts by multiplying the $i$-th factor by $J$
(here we regard $M$ as the 0-th factor).
Thus we have a natural action of $W$ on 
$F_{n,p,\mu}(M)$. 

Define elements $\tilde{y}_{k}\in \End(F_{n,p,\mu}(M))$ as follows: 
 
\begin{equation}\label{new-y}
\tilde{y}_{k} =-\sum_{i|j}E_{ij}\otimes(E_{ji})_{k}, 
\text{ for } k=1,\ldots, n,
\end{equation}
where the first component acts 
on $M$ and the second component acts on the $k$-th factor
of the tensor product.

The main result of this section is the following theorem. 

\begin{theorem}\label{affi} The above action of $W$ and 
the elements $\tilde{y}_k$ given by \eqref{new-y} combine into 
a representation of the degenerate affine Hecke algebra 
$\mathcal{H}(\kappa_1,\kappa_2)$ (in the presentation of 
Lemma \ref{rel-tildey1}) on the space $F_{n,p,\mu}(M)$,
with
\begin{equation}\label{par-rel-1}
\kappa_1=1,\ \kappa_2=p-q-\mu N.
\end{equation}
So we have an functor $F_{n,p,\mu}$ 
from the the category of 
$\mathfrak{g}$-modules to the category of 
representations of type $BC_{n}$ dAHA with such parameters.
If we restrict this functor on the category of Harish-Chandra modules,
we get an exact functor. 
\end{theorem}

\begin{proof} 
Our job is to show that 
the elements $\tilde{y}_{k}$, $S_{k}$ and $\gamma_{n}$ 
satisfy the relations in Lemma \ref{rel-tildey1}.
We only need to prove the commutation 
relation between the elements $\tilde{y}_{k}$, since the other
relations are trivial.

Let $a\neq b$ and $\delta_{ij}$ be the identity matrix, 
then we have
\begin{eqnarray*}%brfu
&&[\tilde{y}_a, \tilde{y}_b]\\
&=& \sum_{il|j}E_{il}\otimes( E_{ji})_a\otimes( E_{lj})_b
 - \sum_{i|kj}E_{kj}\otimes( E_{ji})_a\otimes( E_{ik})_b\\
&=&\sum_{il|j}(E_{il}-\frac{I_N}{N}\delta_{il})\otimes(E_{ji})_{a}\otimes(E_{lj})_{b}
-\sum_{i|kj}(E_{kj}-\frac{I_N}{N}\delta_{kj})\otimes(E_{ji})_{a}\otimes(E_{ik})_{b}\\
&&\text{(By the $\mu$-invariance and lemma \ref{chi-mu})}\\
&=&((q-p)+\mu N)(\sum_{i\leq p, j>p}1\otimes(E_{ji})_{a}\otimes(E_{ij})_{b}
-\sum_{i>p, j\leq p}1\otimes(E_{ji})_{a}\otimes(E_{ij})_{b})\\
&&-\sum_{a\neq c\neq b}(\sum_{il|j}1\otimes(E_{ji})_{a}\otimes(E_{lj})_{b}\otimes(E_{il})_{c}
-\sum_{i|kj}1\otimes(E_{ji})_{a}\otimes(E_{ik})_{b}\otimes(E_{kj})_{c})
\end{eqnarray*}
\begin{eqnarray*}
&=&\frac{p-q-\mu N}{2}S_{ab}(\gamma_{b}-\gamma_{a})
-\frac{1}{4}\sum_{a\neq c\neq b}(1-\gamma_{a}\gamma_{b}-\gamma_{a}\gamma_{c}+\gamma_{b}\gamma_{c})S_{ac}S_{bc}\\
&&\qquad+\frac{1}{4}\sum_{a\neq c\neq b}(1-\gamma_{a}\gamma_{b}+\gamma_{a}\gamma_{c}-\gamma_{b}\gamma_{c})S_{bc}S_{ac}.
 \end{eqnarray*}

Comparing this to the relation in Lemma \ref{rel-tildey1}, we get
the result.
\end{proof}

\subsection{Example}

Consider the example $p=q=1$, $N=2$. 
Thus, $\kappa_1=1$, $\kappa_2=-\mu N$. 
For the module $M$, let us take the module 
${\mathcal F}_{\lambda,\nu}$ of tensor fields 
$p(z)z^{\nu/2} (dz/z)^\lambda$, where $p$ is a Laurent polynomial; 
the Lie algebra $\mathfrak{gl}_2$ acts in it by infinitesimal fractional
linear transformations of $z$. Then we get 
$F_{n,p,\mu}({\mathcal F}_{\lambda,\mu-n})=Y_\lambda$, a representation 
of $\mathcal{H}(\kappa_1,\kappa_2)$ of dimension $2^n$, which is isomorphic to
$V^{\otimes n}$ as a $W$-module. The structure of $Y_\lambda$ as
a dAHA-module is discussed in the recent paper \cite{M}. 

\section{Construction of the functor $F^\lambda_{n,p,\mu}$}
%%%%%%%%%%%%%%%%%%%%%%%%%%%%%%%%
%%%%%%%%%%%%%%%%%%%%%%%%%%%%%%%%

\subsection{The main theorem}

Let $\lambda\in \Bbb C$. 
For $x\in \mathfrak{g}$, let $L_x$ denote the vector field on 
$G$ generated by the left action of $x$.
Thus, $(L_xf)(A)=\frac{d}{dt}|_{t=0}f(e^{tx}A)$ for a function
$f$. Note that $L_{[x,y]}=-[L_x,L_y]$ (the minus sign comes from the
fact that left multiplication by elements of $G$ gives rise to a
{\it right} action of $G$ on functions on $G$). 

Let $\mathcal{D}^{\lambda}(G/K)$ be the sheaf of differential
operators on $G/K$, twisted by the character $\lambda\chi$. 
Local sections of ${\mathcal D}^{\lambda}(G/K)$ act
naturally on $\lambda\chi$-twisted functions on $G/K$, i.e. 
analytic functions $f$ on a small open set $U\subset G$ 
such that \linebreak $R_zf=\lambda\chi(z)f$, $z\in \kk$, where $R_z$ is the
left invariant vector field corresponding to the right
translation by $z$. This action is faithful. Note that we can regard elements
$L_x$ as global sections of  $\mathcal{D}^{\lambda}(G/K)$,
with the same commutation law $[L_x,L_y]=-L_{[x,y]}$. 

Let $M$ be a $\mathcal{D}^{\lambda}(G/K)$-module. 
Then $M$ is naturally a $\mathfrak{g}$-module, via the vector
fields $L_x$. Define
 $$
F^\lambda_{n,p,\mu}(M)=
(M\otimes V^{\otimes n})^{\kk_0,\mu}.
$$ 
Then $F^\lambda_{n,p,\mu}(M)$ is a $W$-module as in the previous
section. 

For $k=1,\ldots, n$, define the following linear operators
on the space $F^\lambda_{n,p,\mu}(M)$:
\begin{eqnarray}\label{rep-X-y}
X_{k}=\sum_{i,j}(AJA^{-1}J)_{ij}\otimes(E_{ij})_{k},\quad
\tilde{y}_{k}=\sum_{i|j}L_{ij}\otimes(E_{ji})_{k},
\end{eqnarray}
where $(AJA^{-1}J)_{ij}$ is the function of $A\in G/K$ which takes the $ij$
-th element of $AJA^{-1}J$, 
$L_{ij}=L_{E_{ij}}$, and the second component acts on the 
$k$-th factor in $V^{\otimes n}$. 
 From now on, we write $X=AJA^{-1}J$ and $X^{-1}=JAJA^{-1}$.
Thus we have $JX=X^{-1}J$. 

The main result of this section is the following theorem. 

\begin{theorem}\label{daffi} The above action of $W$ and 
the elements in \eqref{rep-X-y} combine into
a representation of the degenerate double affine Hecke algebra 
$\HH(t,k_{1},k_{2},k_{3})$ (in the presentation of Lemma \ref{rel-tildey})
on the space $F^\lambda_{n,p,\mu}(M)$, with  
\begin{equation}\label{par-rel-2}
t=\dfrac{2n}{N}+(\lambda+\mu)(q-p),\quad k_{1}=1,\quad k_{2}=p-q-\lambda N,\quad k_{3}=(\lambda-\mu)N.
\end{equation}
So we have an functor $F^\lambda_{n,p,\mu}$ 
from the the category of  
$\mathcal{D}^{\lambda}(G/K)$-modules 
to the category of representations of the type $BC_{n}$ 
dDAHA with such parameters.
If we restrict this functor on the category of $K$-monodromic twisted
$\mathcal{D}$-modules, we get an exact functor.
\end{theorem}

Note that the restriction of the representation 
$F^\lambda_{n,p,\mu}(M)$ to the affine Hecke algebra
$\mathcal{H}$ clearly
coincides with the representation of Theorem \ref{affi}. 

\vskip .05in
{\bf Remark.} This construction is parallel to the
similar construction in type $A$ performed in \cite{CEE}. 
The most important difference is that here we use 
the functions $(AJA^{-1}J)_{ij}$ on $G/K$ instead of the functions 
$A_{ij}$ on $G$ used in \cite{CEE}. The motivation for this is that 
the matrix elements of $AJA^{-1}J$ are ``the simplest''
nonconstant algebraic functions on $G/K$, similarly to 
how the matrix elements of $A$ are ``the simplest'' nonconstant
algebraic functions on $G$.  
\vskip .05in

The rest of this section is devoted to the proof of
Theorem \ref{daffi}. 

\subsection{Proof of Theorem \ref{daffi}}

Our job is to show that 
the elements $X_{k}$, 
$\tilde{y}_{k}$, $S_{k}$ and $\gamma_{n}$ satisfy the relations in Lemma \ref{rel-tildey}.

First of all, the relations in Lemma \ref{rel-tildey} which
don't involve $X_i$ can be established as in the proof of Theorem \ref{affi}
(as \eqref{par-rel-1} is compatible with \eqref{par-rel-2}). 

Second, there are some trivial relations:
\begin{eqnarray*}
&[X_{i},X_{j}]=0, \qquad [\gamma_{i}, X_{j}]=0,\quad(i\neq j),\\
&[S_{i},X_{j}]=0,\quad(j\neq i, i+1),\qquad S_{i}X_{i}-X_{i+1}S_{i}=0,\\
&[\gamma_{i},\tilde{y}_{j}]=0,\quad(j\neq i), \qquad \gamma_{i}\tilde{y}_{i}+\tilde{y}_{i}\gamma_{i}=0,
 \end{eqnarray*} 
and since $JX=X^{-1}J$, we have 
$$
\gamma_{i}X_{i}=X_{i}^{-1}\gamma_{i}.
$$
 
Third, we have the following result.
\begin{lemma}\label{rel-Xy}
We have the following commutation relations: if $m\neq k$ then
\begin{eqnarray*}
&[\tilde{y}_{m},X_{k}] & =  \frac{1}{2}(X_{k}+X_{m})S_{mk}-\frac{1}{2}(X_{k}+X_{m}^{-1})S_{mk}\gamma_{m}\gamma_{k},\\
&[\tilde{y}_{m},X_{k}^{-1}] & =  -\frac{1}{2}(X_{k}^{-1}+X_{m}^{-1})S_{mk}+\frac{1}{2}(X_{k}^{-1}+X_{m})S_{mk}\gamma_{m}\gamma_{k}.
\end{eqnarray*}
So we have 
\begin{equation*}
[\tilde{y}_{m},X_{k}+X_{k}^{-1}] =  \frac{1}{2}(X_{k}-X_{k}^{-1}+X_{m}-X_{m}^{-1})S_{mk}+\frac{1}{2}(X_{k}^{-1}-X_{k}+X_{m}-X_{m}^{-1})S_{mk}\gamma_{m}\gamma_{k}.
\end{equation*}
\end{lemma}

\begin{proof}
The proof is by direct computation.
First, we have for $r\leq p<s$ or $s\leq p<r$

\begin{eqnarray*}
L_{rs}(X)_{ij}
=\delta_{sj}(X)_{ir}+\delta_{ri}(X)_{sj}.
\end{eqnarray*}
So
\begin{eqnarray*}
&&[\tilde{y}_{m},X_{k}]\\
&=& \sum_{r|s}\sum_{i,j}L_{rs}(X)_{ij}\otimes(E_{sr})_{m}\otimes(E_{ij})_{k}\\
&=&\sum_{r|s}\sum_{i}(X)_{ir}\otimes(E_{sr})_{m}\otimes(E_{is})_{k}
+\sum_{r|s}\sum_{j}(X)_{sj}\otimes(E_{sr})_{m}\otimes(E_{rj})_{k}\\
&= &\frac{1}{2}(X_{k}+X_{m})S_{mk}-\frac{1}{2}(X_{k}+X_{m}^{-1})S_{mk}\gamma_{m}\gamma_{k}.
\end{eqnarray*}
By a similar method, we can get the other identities.
\end{proof}

Thus, we only need to show that $X_{m}$ and $\tilde{y}_{m}$
satisfy v) in Lemma \ref{rel-tildey} if the parameters satisfy
\eqref{par-rel-2}.  Instead of computing $[\tilde{y}_{m},X_{m}]$,
we will compute \linebreak $[\tilde{y}_{m},X_{m}+X_{m}^{-1}]$ and
$[\tilde{y}_{m},X_{m}-X_{m}^{-1}]$.

\subsubsection{Computing $[\tilde{y}_{m},X_{m}+X_{m}^{-1}]$}

Let us define 
$$
T=\tr(X)=\sum_{i}(X)_{ii}\otimes 1.
$$

Suppose $X=\left(\begin{array}{cc}A_1 & A_2 \\A_3 &
A_4\end{array}\right)$ 
where $A_{1}$ is a $p$ by $p$ matrix.
Then 
\begin{equation}
\tr(A_{1})=\tr(X(J+1)/2),\qquad \tr(A_{4})=\tr(X(1-J)/2).
\end{equation}
But $\tr(XJ)=\tr(AJA^{-1})=\tr(J)=p-q$, so we get 
\begin{equation}\label{tr12}
\tr(A_{1})=\frac{T+p-q}{2},\qquad \tr(A_{4})=\frac{T-p+q}{2}.
\end{equation}

\begin{lemma}\label{lemma-T}
We have the relation
$$\sum_{m}(X_{m}+X_{m}^{-1})=(\frac{2n}{N}+(\lambda+\mu)(q-p))T+(\lambda+\mu)(p^{2}-q^{2}).$$
\end{lemma}

\begin{proof}
Since $(X)_{ij}=-(X^{-1})_{ij}$ unless $i,j\leq p$ or $i,j>p$, 
and $(X)_{ij}=(X^{-1})_{ij}$ if $i,j\leq p$ or $i,j>p$,
we have 
$$
X_{m}+X_{m}^{-1}=\sum_{ij}(X+X^{-1})_{ij}\otimes (E_{ij})_{m}.
$$

Thus we have 
\begin{eqnarray}\label{sumx}
&& \sum_m(X_{m}+X_{m}^{-1})\\\nonumber
& = & \sum_m\sum_{ij}(X+X^{-1})_{ij}
\otimes (E_{ij}-\frac{I_N}{N}\delta_{ij})_{m}+\sum_m\sum_{i}
(X+X^{-1})_{ii}\otimes (\frac{I_N}{N})_{m}\\\nonumber
&&\text{(By the $\mu$-invariance and Lemma \ref{chi-mu})}\\\nonumber
& = &Y +\frac{2n}{N}T+\sum_{i\leq p}\mu q(X+X^{-1})_{ii}
\otimes 1-\sum_{i> p}\mu 
p(X+X^{-1})_{ii}\otimes 1\\\nonumber
&&\text{(By (\ref{tr12}))}\\\nonumber
&=&Y+(\frac{2n}{N}+\mu(q-p))T+\mu(p^{2}-q^{2}),
\end{eqnarray}
where
$Y=\sum_{ij}(X+X^{-1})_{ij}L_{E_{ij}}\otimes
1$. 

It remains to calculate the expression $Y$ 
in the algebra ${\mathcal D}^{\lambda}(G/K)$.
We can calculate $Y$ by acting with it on
$\lambda\chi$-twisted functions $f$ on $G/K$. 

We have 
\begin{eqnarray*}
(Yf)(A)&=&\frac{d}{dt}|_{t=0}f(A+t(X+X^{-1})A)\\
&=&
\frac{d}{dt}|_{t=0}f(A+tAJA^{-1}JA+tJAJ)\\
&=&
\frac{d}{dt}|_{t=0}f(A+tA(JA^{-1}JA+A^{-1}JAJ))\\
&=&
\frac{d}{dt}|_{t=0}f(A+tA(X_*+X_*^{-1})),
\end{eqnarray*}
where $X_*:=JA^{-1}JA$. Now, $X_*+X_*^{-1}\in \kk$,
so we have 
$$
Yf=\lambda\chi(X_*+X_*^{-1})f=\lambda((q-p)T+(p^2-q^2))f.
$$
Combining this with formula (\ref{sumx}), we obtain the statement
of the lemma. 
\end{proof}

Notice that
\begin{eqnarray*}
[\tilde{y}_{m}, T]
& = & \sum_{r|s}L_{rs}(\tr(X))\otimes (E_{sr})_{m}\\
& = &\sum_{r|s}((X)_{sr}-(X^{-1})_{sr})\otimes (E_{sr})_{m}\\
& = & X_{m}-X_{m}^{-1}.
\end{eqnarray*}

So from Lemma \ref{rel-Xy} and Lemma \ref{lemma-T}, we have
\begin{eqnarray*}
[\tilde{y}_{m},X_{m}+X_{m}^{-1}]=-\sum_{k\neq
m}[\tilde{y}_{m},X_{k}+X_{k}^{-1}]+(\frac{2n}{N}+
(\lambda+\mu)(q-p))[\tilde{y}_{m},T].
\end{eqnarray*}
Thus, we have obtained

\begin{lemma}\label{xplusx}
\begin{eqnarray}\label{rel-y-x+x}
&&[\tilde{y}_{m},X_{m}+X_{m}^{-1}]\\\nonumber
& = &(\frac{2n}{N}+(\lambda+\mu)(q-p))(X_{m}-X_{m}^{-1})-
\frac{1}{2}\sum_{k\neq m}(X_{k}-X_{k}^{-1}+X_{m}-X_{m}^{-1})S_{mk}\\\nonumber
&&\qquad+\frac{1}{2}\sum_{k\neq m}(X_{k}-X_{k}^{-1}-X_{m}+X_{m}^{-1})S_{mk}\gamma_{m}\gamma_{k}.
\end{eqnarray}
\end{lemma}

\subsubsection{Computing $[\tilde{y}_{m},X_{m}-X_{m}^{-1}]$}

At first, we need the following lemmas for the future computation.

\begin{lemma}\label{lemma-qpsum}
We have the equality
\begin{eqnarray*}
&&q\sum_{s\leq p}\sum_{j}(X+X^{-1})_{sj}\otimes (E_{sj})_{m}+  
p\sum_{s>p}\sum_{j}(X+X^{-1})_{sj}\otimes (E_{sj})_{m}\\
&=&\frac{1}{2}(N+(q-p)\gamma_{m})(X_{m}+X_{m}^{-1}).
\end{eqnarray*}
\end{lemma}
\begin{proof}
By direct computation.
\end{proof}

\begin{lemma}\label{lemma-vf-equal}
In $\mathcal{D}^{\lambda}(G/K)$, we have for $r,j\leq p$ or $r,j>p$,
$$L_{[X-X^{-1},E_{rj}]}=-L_{\{X+X^{-1},E_{rj}\}}+2\lambda\chi(Q_{rj}),$$
where $\{a,b\}=ab+ba$ and $Q_{rj}=JA^{-1}JE_{rj}A+A^{-1}E_{rj}JAJ$.
\end{lemma}
\begin{proof}
Let $f(A)$ be a $\lambda\chi$-twisted function on $G/K$,
i.e. an analytic function on a small open set $U\subset G$ such
that $R_zf=\lambda\chi(z)f$, $z\in \kk$.  
Then we have 
\begin{eqnarray*}
&&L_{[X-X^{-1},E_{rj}]}f(A)\\
&=&\frac{d}{dt}|_{t=0}f(A+t(XE_{rj}+E_{rj}X^{-1})A
-t(X^{-1}E_{rj}+E_{rj}X)A).
\end{eqnarray*}
Notice that
\begin{eqnarray*}
&&f(A+t(XE_{rj}+E_{rj}X^{-1})A
-t(X^{-1}E_{rj}+E_{rj}X)A)\\
&=&f(A+2tAQ_{rj}-t(XE_{rj}+E_{rj}X^{-1})A
-t(X^{-1}E_{rj}+E_{rj}X)A),
\end{eqnarray*}
and $Q_{rj}$ is an element of $\kk$.

So we have
\begin{eqnarray*}
&&\frac{d}{dt}|_{t=0}f(A+tAQ_{rj}-t(X^{-1}E_{rj}+E_{rj}X)A)\\
&=&\frac{d}{dt}|_{t=0}f(A-t(X^{-1}E_{rj}+E_{rj}X)A-t(XE_{rj}+E_{rj}X^{-1})A)+\frac{d}{dt}|_{t=0}f(A+2tAQ_{rj})\\
&=&-L_{\{X+X^{-1},E_{rj}\}}f(A)+2\lambda\chi(Q_{rj})f(A).
\end{eqnarray*}
Thus we get the lemma.
\end{proof}

Now let us compute $[\tilde{y}_{m},X_{m}-X_{m}^{-1}]$.
By the definition and Lemma \ref{lemma-vf-equal}, we have
\begin{eqnarray}\label{eqn-1}
&&[\tilde{y}_{m},X_{m}-X_{m}^{-1}]\\\nonumber
& = &  \sum_{r|s}\sum_{i,j}L_{rs}((X)_{ij}-(X^{-1})_{ij})\otimes (E_{sr}E_{ij})_{m}
-\sum_{rj}L_{\{X+X^{-1},E_{jr}\}}\otimes (E_{rj})_{m}\\
&&\qquad+
2\lambda\sum_{rj}\chi(JA^{-1}JE_{rj}A+A^{-1}E_{rj}JAJ)
\otimes (E_{jr})_{m}.\nonumber
\end{eqnarray}

Since we have
\begin{eqnarray*}
L_{rs}((X)_{ij}-(X^{-1})_{ij})
= (X)_{ir}\delta_{sj}+(X)_{sj}\delta_{ir}
          +(X^{-1})_{ir}\delta_{sj}+(X^{-1})_{sj}\delta_{ir},
\end{eqnarray*}
by Lemma \ref{lemma-qpsum}, the first summand of \eqref{eqn-1} is
\begin{eqnarray*}
\frac{1}{2}(N+(q-p)\gamma_{m})(X_{m}+X_{m}^{-1})+(1+\gamma_{m})\frac{T-p+q}{2}+(1-\gamma_{m})\frac{T+p-q}{2}.   
\end{eqnarray*}

Now let us compute the second summand of \eqref{eqn-1}.
By definition, we have
\begin{eqnarray*}
&&-\sum_{rj}L_{\{X+X^{-1},E_{rj}\}}\otimes (E_{jr})_{m}\\\nonumber
&=&-\sum_{ijr}\left((X+X^{-1})_{ir}L_{E_{ij}}+(X+X^{-1})_{ji}L_{E_{ri}}\right)\otimes(E_{jr})_{m}\\\nonumber
&=&-\sum_{ijr}\left((X+X^{-1})_{ir}L_{E_{ij}-\frac{I_N}{N}\delta_{ij}}+(X+X^{-1})_{ji}L_{E_{ri}-\frac{I_N}{N}\delta_{ri}}\right)\otimes(E_{jr})_{m}\\
&=&-\sum_{ijr}(X+X^{-1})_{ir}\otimes(E_{jr}E_{ij})_{m}-\sum_{k\neq m}\sum_{ijr}(X+X^{-1})_{ir}\otimes(E_{ij})_{k}\otimes(E_{jr})_{m}\\
&&-\sum_{ijr}(X+X^{-1})_{ji}\otimes(E_{jr}E_{ri})_{m}-\sum_{k\neq m}\sum_{ijr}(X+X^{-1})_{ji}\otimes(E_{ri})_{k}\otimes(E_{jr})_{m}\\
&&+\sum_{ijr}(X+X^{-1})_{ir}\otimes(E_{jr}\frac{I_N}{N}\delta_{ij})_{m}+\sum_{k\neq m}\sum_{ijr}(X+X^{-1})_{ir}\otimes(\frac{I_N}{N}\delta_{ij})_{k}\otimes(E_{jr})_{m}\\
&&+\sum_{ijr}(X+X^{-1})_{ji}\otimes(E_{jr}\frac{I_N}{N}\delta_{ri})_{m}+\sum_{k\neq m}\sum_{ijr}(X+X^{-1})_{ji}\otimes(\frac{I_N}{N}\delta_{ri})_{k}\otimes(E_{jr})_{m}\\
&&+\mu((q-p)+N\gamma_{m})(X_{m}+X_{m}^{-1})
\end{eqnarray*}
\begin{eqnarray*}
&=&-T+(q-p)\gamma_{m}-\frac{1}{2}(X_{m}+X_{m}^{-1})(N+(p-q)\gamma_{m})\\
&&\quad-\frac{1}{2}\sum_{k\neq m}(X_{m}+X_{m}^{-1}+X_{k}+X_{k}^{-1})S_{km}(1+\gamma_{k}\gamma_{m})\\
&&\qquad+\frac{2n(X_{m}+X_{m}^{-1})}{N}+\mu((q-p)+N\gamma_{m})(X_{m}+X_{m}^{-1}).
\end{eqnarray*}

Now let us compute the third summand of \eqref{eqn-1}.

\begin{lemma}\label{chiS}
For $r,j\leq p$ or $r,j>p$, we have:
\begin{eqnarray*}
&&\lambda\chi(JA^{-1}JE_{rj}A+A^{-1}E_{rj}JAJ)\\
&=&\left\{\begin{array}{cc}
\dfrac{\lambda(q-p)}{2}(X+X^{-1})_{jr}+\lambda N\delta_{rj}, & r,j\leq p; \\
&\\
\dfrac{\lambda(q-p)}{2}(X+X^{-1})_{jr}-\lambda N\delta_{rj}, & r,j> p. \end{array}\right.
\end{eqnarray*} 
\end{lemma}
\begin{proof}
Let us denote 
$$
B=A^{-1}E_{rj}JA=\left(\begin{array}{cc}B_{1} & B_{2} \\B_{3} & B_{4}\end{array}\right), \text{ where $B_{1}$ is a $p$ by $p$ matrix.}
$$
Then 
$$
JA^{-1}JE_{rj}A+A^{-1}E_{rj}JAJ
=JB+BJ
=\left(\begin{array}{cc}2B_{1} & \\ & -2B_{4}\end{array}\right).
$$

Then we have
$$
2\tr(B_{1})-2\tr(B_{4})
=(X+X^{-1})_{jr}.
$$

On the other hand we have
$$
\tr(B)=\tr(B_{1})+\tr(B_{4})=\tr(A^{-1}E_{rj}JA)=
\left\{\begin{array}{cc}
\delta_{rj}, & rj\leq p; \\
 -\delta_{rj},& rj>p. \end{array}\right.
$$

Then we have for $j,r\leq p$,
$$\tr(B_{1})=\frac{1}{4}(X+X^{-1})_{jr}
+\frac{1}{2}\delta_{rj},\quad \tr(B_{4})=-\frac{1}{4}(X+X^{-1})_{jr}
+\frac{1}{2}\delta_{rj},$$
for $j,r> p$,
$$\tr(B_{1})=\frac{1}{4}(X+X^{-1})_{jr}
-\frac{1}{2}\delta_{rj},\quad \tr(B_{4})=-\frac{1}{4}(X+X^{-1})_{jr}
-\frac{1}{2}\delta_{rj}.$$

So we get the lemma.
\end{proof}

From Lemma \ref{chiS}, we have:
\begin{eqnarray*}
&&\sum_{rj}\chi(JA^{-1}JE_{rj}A+A^{-1}E_{rj}JAJ)\otimes (E_{jr})_{m}\\
&&\qquad
=\frac{1}{2}(q-p)(X_{m}+X_{m}^{-1})+N\gamma_{m}.
\end{eqnarray*}

Thus, combining the above formulas, we have

\begin{lemma}\label{xminx}
\begin{eqnarray}\label{rel-y-x-x}
&&[\tilde{y}_{m},X_{m}-X_{m}^{-1}]\\\nonumber
&=&(\frac{2n}{N}+(\lambda+\mu)(q-p))(X_{m}+X_{m}^{-1})-\frac{1}{2}\sum_{k\neq m}(X_{k}+X_{k}^{-1})(1+\gamma_{m}\gamma_{k})S_{km}\\\nonumber
&&\qquad-\frac{1}{2}\sum_{k\neq m}(X_{m}+X_{m}^{-1})(1+\gamma_{m}\gamma_{k})S_{km}\\\nonumber
&&\qquad+((q-p)+\mu N)
\gamma_{m}(X_{m}+X_{m}^{-1})+2((q-p)+\lambda N)\gamma_{m}.
\end{eqnarray}
\end{lemma}

\subsubsection{Conclusion}

Adding equations (\ref{rel-y-x-x}) and (\ref{rel-y-x+x}),
and comparing with Lemma \ref{rel-tildey}, we conclude the
proof of Theorem \ref{daffi}. 

\section{Action of $F^\lambda_{n,p,\mu}$ on some subcategories}

As we mentioned, the functors $F_{n,p,\mu}$,
$F_{n,p,\mu}^\lambda$ factor through
the category of modules $M$ on which the action of $\kk$ is
locally finite; more precisely,
$F_{n,p,\mu}(M)=F_{n,p,\mu}(M_f)$, 
where $M_f$ is the locally finite part of $M$ under $\kk$. 
 
Now let $M$ be a $\mathcal{D}^{\lambda}(G/K)$-module,
locally finite under $\kk$. 
The support of such a $\mathcal{D}$-module is a $K$-invariant subset of 
$G/K$, i.e. a union of $K$-orbits. Recall that closed $K$-orbits
of $G/K$ are labeled by the points of the categorical quotient 
$K\backslash G/K$, i.e. the spectrum of the ring 
$R_{p,q}={\mathcal O}(G/K)^K$. 
For every point $\psi\in K\backslash G/K$ ($\psi: R_{p,q}\to \Bbb C$), we
can define the subcategory $D^\lambda(\psi)$ of the category 
of ${\mathcal D}^{\lambda}$-modules on $G/K$ which are 
set-theoretically supported on the preimage of $\psi$ in $G/K$,
i.e. those on which $R_{p,q}$ acts with generalized eigenvalue
$\psi$. 

On the other hand, let $\Bbb T=\Bbb C^n/\Bbb Z^n=\h/Q^\vee$, and 
let $\beta\in \Bbb T/W$. Then we can define the category 
${\mathcal O}_\beta$ of modules over the dDAHA $\HH$ 
on which the subalgebra $\Bbb C[P]^W\subset \HH$ acts with 
generalized eigenvalue $\beta$. 

The following theorem tells us how the functor
$F^\lambda_{n,p,\mu}$ relates $\psi$ and $\beta$. 

\begin{theorem}\label{thet}
The functor $F^\lambda_{n,p,\mu}$ maps 
$D^{\lambda}(\psi)$ to ${\mathcal O}_\beta$, where $\beta=\theta(\psi)$,\linebreak
and $\theta: K\backslash G/K\to\Bbb T/W$ is the regular map defined by the
formula 
$$
\theta^*(\sum_{m=1}^n g(X_m))=ng(1)+
(\frac{n}{N}+\frac{1}{2}(\lambda+\mu)(q-p))\tr(g(X)-g(1)),
$$
where $g$ is a 
Laurent polynomial in one variable such that $g(Z)=g(Z^{-1})$.
\end{theorem} 

\begin{proof}
The proof is obtained by generalizing the proof of 
Lemma \ref{lemma-T}.
We'll need the following lemma. 

\begin{lemma} 
$\tr(X^sJ)=p-q$ for any $s\in \Bbb Z$. 
\end{lemma}
\begin{proof}
It's easy to see that $X^sJ$ is conjugate to $J$. 
\end{proof}

The lemma implies that 
$$
\tr(X^s(J+1)/2)=\frac{\tr(X^s)+p-q}{2},
$$ 
and 
$$
\tr(X^s(1-J)/2)=\frac{\tr(X^s)-p+q}{2}.
$$ 

Thus we have 
\begin{eqnarray}\label{sumx1}
&& \sum_m g(X_{m})\\\nonumber
& = & \sum_m\sum_{ij}g(X)_{ij}
\otimes (E_{ij}-\frac{I_N}{N}\delta_{ij})_{m}+\sum_m\sum_{i}
            g(X)_{ii}\otimes (\frac{I_N}{N})_{m}\\\nonumber
&&\text{(By the $\mu$-invariance and Lemma \ref{chi-mu})}\\\nonumber
& = &Y_g +\frac{n}{N}\tr(g(X))+
\sum_{i\leq p}\mu qg(X)_{ii}\otimes 1-\sum_{i> p}\mu pg(X)_{ii}\otimes 1\\\nonumber
&=&Y_g+(\frac{n}{N}+\frac{1}{2}\mu(q-p))\tr(g(X))+\frac{1}{2}
\mu(p^{2}-q^{2})g(1),
\end{eqnarray}
where
$Y_g=\sum_{ij}g(X)_{ij}L_{E_{ij}}\otimes
1$. 

It remains to calculate the expression $Y_g$ 
in the algebra ${\mathcal D}^{\lambda}(G/K)$.
We can calculate $Y_g$ by acting with it on
$\lambda\chi$-twisted functions $f$ on $G/K$. 

We have 
\begin{eqnarray*}
(Y_gf)(A)&=&\frac{d}{dt}|_{t=0}f(A+tg(X)A)\\
&=&\frac{d}{dt}|_{t=0}f(A+tAg(X_*)),
\end{eqnarray*}
where $X_*:=JA^{-1}JA$. Now, $X_*+X_*^{-1}\in \kk$,
so we have 
$$
Y_gf=\lambda\chi(g(X_*))f=
\frac{1}{2}\lambda((q-p)\tr(g(X))+\frac{1}{2}(p^2-q^2)g(1))f.
$$
Combining this with formula (\ref{sumx1}), we obtain the statement
of the theorem. 
\end{proof}

\begin{remark} 
In particular, Theorem \ref{thet} implies that 
$\theta(1)=1$, where $1\in K\backslash G/K$ is the double coset
of $1$, and $1\in \Bbb T/W$ is the image of the unit of the group $\Bbb T$.  
Thus the functor $F^\lambda_{n,p,\mu}$ maps the category $D^{\lambda}(1)$ 
to the category ${\mathcal O}_1$. Note that $D^{\lambda}(1)$ is the
category of twisted $\mathcal{D}$-modules supported on the ``unipotent
variety'' in $G/K$ (which is equivalent to the category 
on $\mathcal{D}$-modules 
on $\mathfrak{g}/\kk$ supported on the nilpotent cone), 
and ${\mathcal O}_1$ is the category of $\HH$-modules on which
$X_i$ act unipotently (which is equivalent to category ${\mathcal
O}$ for the rational Cherednik algebra of type $B_n$). 
\end{remark}

\begin{remark}
Another very interesting question is how the functor $F_{n,p,\mu}$
transforms the central characters, i.e. how the central character
of a Harish-Chandra module $M$ is related to the central
character of the dAHA module $F_{n,p,\mu}(M)$. This question is
discussed in the paper \cite{M} (for $n=1$). 
\end{remark} 

\section*{acknowledgments}
%%%%%%%%%%%%%%%%%%%%

{The work of the first author was  partially supported by the NSF grant 
DMS-0504847. The work of the second and the third author was supported by
the Summer Program of Undergraduate Research in the Department of
Mathematics at MIT. We thank Ju-Lee Kim, David Vogan, and Ting Xue for useful discussions.}

%%%%%%%%%%%%%%%%%%%%

\end{document}